\documentclass[a4paper,11pt,reqno]{amsart}
\usepackage[latin1]{inputenc}
\usepackage[english]{babel}
\usepackage{amssymb,amsfonts, amsmath, color}
\usepackage[margin=2.7cm]{geometry}
\usepackage[dvipsnames]{xcolor}
\usepackage{graphicx, color, enumerate}
\usepackage[latin1]{inputenc}
\usepackage[active]{srcltx}
\usepackage{tikz}
\usepackage{pgf}
\usepackage{etex}
\usepackage{verbatim}
\usepackage{pgfkeys}
\usepackage{amsmath}
\usepackage{amsfonts}
\usepackage{amssymb}
\usepackage{amsthm}
\usepackage{xcolor}
\usepackage{mathrsfs}
\usepackage{geometry}
\usepackage{fancyhdr}
\usepackage[colorlinks, citecolor=blue, linkcolor=red]{hyperref}

\newtheorem{theorem}{Theorem}[section]
\newtheorem{proposition}[theorem]{Proposition}

\newtheorem{remark}[theorem]{Remark}
\newtheorem{definition}[theorem]{Definition}

\newcommand{\bcl}{\begin{center}}
\newcommand{\ecl}{\end{center}}
\newcommand{\brl}{\begin{right}}
\newcommand{\erl}{\end{right}}
\newcommand{\ben}{\begin{enumerate}}
\newcommand{\een}{\end{enumerate}}
\newcommand{\overliner}{\begin{array}}
\newcommand{\earr}{\end{array}}
\newcommand{\btab}{\begin{tabular}}
\newcommand{\etab}{\end{tabular}}
\newcommand{\bdoc}{\begin{document}}
\newcommand{\edoc}{\end{document}}
\newcommand{\beqy}{\begin{eqnarray}}
\newcommand{\eeqy}{\end{eqnarray}}

\newcommand{\beqi}{\begin{eqnarray*}}
\newcommand{\eeqi}{\end{eqnarray*}}
\newcommand{\bitem}{\begin{itemize}}

\newcommand{\eitem}{\end{itemize}}
\newcommand{\nln}{\newline}
\newcommand{\newt}{\newtheorem}
\newcommand{\pa}{\partial}
\newcommand{\re}{{I\!\!R}}
\newcommand{\Rn}{\R^N}
\newcommand{\xr}{x\in\R }
\newcommand{\x}{\times}
\newcommand{\dyle}{\displaystyle}
\newcommand{\ene}{{I\!\!N}}
\newcommand{\irn}{\int\limits_{\R^N}}
\newcommand{\io}{\int\limits_{\O}}
\newcommand{\meas}{{\rm meas\,}}
\newcommand{\dif}{\nabla_{xy}}
\newcommand{\sign}{{\rm sign}}
\newcommand{\map}{\longrightarrow }
\newcommand{\imp}{\Longrightarrow }
\renewcommand{\div}{\nabla\cdot }
\newcommand{\sen}{{\rm sen\,}}
\newcommand{\tg}{{\rm tg\,}}
\newcommand{\arcsen}{{\rm arcsen\,}}
\newcommand{\arctg}{{\rm arctg\,}}
\newcommand{\supp}{{\textsl supp\ }}
\newcommand{\ity}{\int_{-\iy}^{+\iy}}
\newcommand{\limit}{\lim\limits}
\newcommand{\limi}{\limit_{n\to\infty}}
\newcommand{\sumi}{\sum\limits_{n=1}^{\infty}}
\newcommand{\ulu}{\underline u}
\newcommand{\ulw}{\underline w}
\newcommand{\ulz}{\underline z}
\newcommand{\ulv}{\underline v}
\newcommand{\uls}{\underline s}
\newcommand{\olu}{\overline u}
\newcommand{\olv}{\overline v}
\newcommand{\ols}{\overline s}
\newcommand{\ob}{\overline\b}
\newcommand{\ovar}{\overline\var}
\newcommand{\wv}{\widetilde v}
\newcommand{\wu}{\widetilde u}
\newcommand{\ws}{\widetilde s}
\renewcommand{\a }{\alpha }
\renewcommand{\b }{\beta }
\newcommand{\g }{\gamma}
\newcommand{\G }{\Gamma }
\renewcommand{\d }{\delta }

\newcommand{\D }{\Delta }
\newcommand{\e }{\varepsilon }
\newcommand{\z }{\zeta }
\renewcommand{\l }{\lambda }
\renewcommand{\L }{\Lambda }
\newcommand{\m }{\mu }
\newcommand{\n }{\nabla }
\newcommand{\s }{\sigma }
\newcommand{\Sig }{\Sigma }
\renewcommand{\t }{\tau }
\newcommand{\var }{\varphi }
\renewcommand{\o }{\omega }
\renewcommand{\O }{\Omega }
\newcommand{\R}{{\mathbb{R}}}
\newcommand{\bC}{{\bf C}}
\newcommand{\bZ}{{\bf Z}}
\newcommand{\bN}{{\bf N}}
\newcommand{\bQ}{{\bf Q}}
\newcommand{\bK}{{\bf K}}
\newcommand{\bI}{{\bf I}}
\newcommand{\bv}{{\bf v}}
\newcommand{\bV}{{\bf V}}
\newcommand{\LL}{\mathcal{L}}
\newcommand{\N}{\mathbb{N}}
\DeclareMathOperator{\suppo}{supp} \DeclareMathOperator{\di}{div}
\newenvironment{Proof}{\Rmovelastskip\vskip12pt
plus 1pt \noindent\em\rm}{\hfill {\qed \hskip .2cm}}

\makeatletter
\@namedef{subjclassname@2020}{\textup{2020} Mathematics Subject Classification}
\makeatother

\begin{document}

\title[Local-nonlocal Fujita-type results]{Global solutions to semilinear parabolic \\
equations driven by mixed local-nonlocal operators}

\author[S.\,Biagi]{Stefano Biagi}
 \author[F.\,Punzo]{Fabio Punzo}
 \author[E.\,Vecchi]{Eugenio Vecchi}

  \address[S.\,Biagi]{Dipartimento di Matematica
 \newline\indent Politecnico di Milano \newline\indent
 Via Bonardi 9, 20133 Milano, Italy}
 \email{stefano.biagi@polimi.it}

 \address[F.\,Punzo]{Dipartimento di Matematica
 \newline\indent Politecnico di Milano \newline\indent
 Via Bonardi 9, 20133 Milano, Italy}
 \email{fabio.punzo@polimi.it}

 \address[E.\,Vecchi]{Dipartimento di Matematica
 \newline\indent Università di Bologna \newline\indent
 Piazza di Porta San Donato 5, 40122 Bologna, Italy}
 \email{eugenio.vecchi2@unibo.it}

\keywords{Mixed local/nonlocal operator, semilinear parabolic equations, heat kernel estimates, test functions.}

\subjclass[2020]{35A01, 35B44, 35K57, 35K58, 35R11}

\date{\today}

\begin{abstract}
We are concerned with the Cauchy problem for the semilinear parabolic equation driven by the mixed local-nonlocal operator $\LL = -\Delta+(-\Delta)^s$, with a power-like source term.
We show that the so-called Fujita phenomenon holds, and the critical value is exactly the same as for the fractional Laplacian.
\end{abstract}
\maketitle

\section{Introduction}
Let $\mathcal L$ be the mixed {\it local-nonlocal} operator $\LL = -\Delta+(-\Delta)^s,$  where $(-\Delta)^s$ stands for the fractional Laplacian of order $s\in (0,1)$.
We investigate global existence and blow-up of solutions to semilinear parabolic equations driven by $\mathcal L$ of the following type:
\begin{equation} \label{eq:PbCauchy}
 \begin{cases}
 \partial_t u + \LL u = u^p & \text{in $\R^N\times(0,+\infty)$} \\
 u = u_0 & \text{in $\R^N$},
 \end{cases}
\end{equation}
where $p>1$ and $u_0$ is a given nonnegative initial datum.
\smallskip

\noindent\emph{Bibliographical notes: global existence and blow-up}.
Global existence and blow-up of solutions have been largely studied in the literature.
Concerning the purely local case $\mathcal L=-\Delta$,
it has been shown in \cite{F}, and in \cite{H} and \cite{KST} for the critical case, that
\begin{enumerate}[a)]
\item if  $1 < p\le 1+\frac{2}N$, any solution of \eqref{eq:PbCauchy} blows up in finite-time, 
provided that $u_0\not\equiv 0$;
\item if $p>1+\frac{2}N$, then there exists a global in time solution of \eqref{eq:PbCauchy}, 
provided that $u_0$ is sufficiently small.
\end{enumerate}
Such a dichotomy is known as the {\it Fujita phenomenon}. We refer, e.g., 
to \cite{BB}, \cite{DL}, \cite{Levine} and the references therein, for a complete account 
about blow-up and global existence of solutions in the purely local
case $\mathcal L=-\Delta$.

This question has been addressed also on Riemannian manifolds when $\mathcal L$ is the Laplace-Beltrami operator; in this direction, some results can be found e.g. in \cite{BPT, GMPu, Sun1, MaMoPu, Punzo, Pu22, WY2, Zhang}. Furthermore, analogue results have also been established for local quasilinear evolution equations (see e.g. \cite{MaMoPu, MeGP1, MeGP2, MeGP3, MeGP4, MeMoPu, Mitidieri2}).

On the other hand, when $\mathcal L=(-\Delta)^s$ in \cite{Sug} it is shown that if $p\leq 1+ \frac{2s}{N}$, then any solution arising from a nontrivial initial datum $u_0$ blows up in finite time (see also \cite{FGS}). Such
a result has been generalized in \cite{LaisterS} for more general source terms. Moreover, in \cite{IKK} (see also \cite{HKN}), for $p>1+\frac{2s}N,$ global in time solutions are considered, and  the asymptotic behavior of solutions as $t\to +\infty$  has been studied.
\medskip

\noindent\emph{Bibliographical notes: mixed local-nonlocal operators}.
Recently, the study of qualitative properties of solutions to partial differential equations, mainly of elliptic but also of parabolic type, driven by the mixed operator $\mathcal L$ has been attracting much attention (see \cite{Biagi2, Biagi3, BiagiBN, Biagi4, Biagi5, Biagi6, BonfE, BFV, BonfII, DeFMin, GarainKinnunen, GarainKinnunen2, GarainKinnunen3, GarainLindgren}). One of the main reasons for this interest
is that
mixed operators of the form $\mathcal L$ have applications in probability; indeed, they  are related to the superposition of different types of sthocastic processes such as a classical random walk and a L\'evy flight. Furthermore, they are exploited to model various phenomena in sciences, such as
the study of optimal animal foraging strategies, see e.g. \cite{DPLV23, DV21} and references therein.
\medskip

\noindent \emph{Description of our results}.
Along the above-described line of research, in the present paper we deal with nonnegative solutions to problem \eqref{eq:PbCauchy}. The main result of this paper will be given in detail in the forthcoming Theorem \ref{thm:Main};
however, we give here a sketchy outline of this result.
In particular, we show that if $p\leq 1+\frac{2s}N$, then problem \eqref{eq:PbCauchy} does not admit any global solution with $u_0\not\equiv 0$. On the other hand, if $p>1+\frac{2s}N$, then there exists a global in time solution, provided that $u_0$ is small enough. We point out that problem  \eqref{eq:PbCauchy} behaves like the problem with $\mathcal L=(-\Delta)^s$; in other terms, for what concerns existence and nonexistence of global in time solutions the mixed local-nonlocal operator has the same character  as the nonlocal operator $(-\Delta)^s$. The proof of the nonexistence of global solutions is based on a test functions argument and on suitable a priori estimates. Furthermore, the global solution is constructed by an iteration method, which exploits in a crucial way the estimates from above for the heat kernel of $\mathcal L$.

\medskip

\noindent\emph{Plan of the paper}. The paper is organized as follows.
In Section \ref{mb} we fix the notation and recall some preliminary results concerning the fractional Laplacian, the operator $\mathcal L$ and the heat kernel of $\mathcal L$.
In Section \ref{mr} we give the precise definition of solution
to problem \eqref{eq:PbCauchy} 
and we state our main existence/non-existence result, which is then proved in Section \ref{proof}.

\section{Mathematical background} \label{mb}\setcounter{equation}{0}
\noindent\textbf{Notation.} Throughout the paper, we will tacitly
exploit all the notation listed below; we thus refer the Reader to this list
for any non-standard notation encountered.
\begin{itemize}
 \item We denote by $\R^+$ (resp.\,$\R^+_0$) the interval $(0,+\infty)$ (resp.\,$[0,+\infty)$).
 \vspace*{0.1cm}

 \item Given any $x_0\in\R^N$ and any $r > 0$, we denote by
 $B_r(x_0)$ the open (Euclidean) ball with centre $x_0$ and radius $r$;
 in the particular case when $x_0 = 0$, we simply write $B_r$.
 \vspace*{0.1cm}

 \item Given any $0<T\leq+\infty$, we denote by $S_T$ the (infinite) strip $\R^N\times (0,T)$;
 in the particular case when $T = +\infty$, we simply write $S$
 in place of $S_{+\infty}$.
 \vspace*{0.1cm}

 \item If $A$ is an arbitrary set in some Euclidean space $\R^m$ (with $m\geq 1$), we denote by $\mathbf{1}_A$ the
 usual indicator function
 of $A$, that is,
 $$\mathbf{1}_A(z) = \begin{cases}
  1& \text{if $z\in A$} \\
  0 & \text{if $z\notin A$}.
  \end{cases}
  $$

 \item We denote by $\mathcal{T}_0$ the set (vector space) of the functions
 $\varphi\in C^\infty(\overline{S})$ for which there exist
 numbers
 $r,T > 0$ (possibly depending on $\varphi$) such that
 $$\text{$\varphi\equiv 0$ on $(\R^N\setminus B_r)\times[T,+\infty)$}.$$
 \item Given any $s\in (0,1)$, we denote by $L_s$ the \emph{tail space}
 \begin{equation*}
 {L}_s(\R^N) :=
 \Big\{f:\R^N\to\R:\,\|u\|_{1,s} := \int_{\R^N}\frac{|f(x)|}{1+|x|^{N+2s}}\,dx<+\infty\Big\}.
\end{equation*}
\item Given any open interval $I\subseteq\R$, any Banach space $(X,\|\cdot\|_X)$
and any $1\leq\theta\leq\infty$, we de\-note by
$L^\theta(I;X)$ the space of the $L^\theta$-functions taking values in $X$, that is,
$$L^\theta(I;X) = \big\{f:I\to X:\,\mathfrak{n}_X(f)(t) := \|f(t)\|_X\in L^\theta(I)\big\}.$$
If $f\in L^\theta(I;X)$, we define $\|f\|_{\theta,I,X} := \|\mathfrak{n}_X(f)\|_{L^\theta(I)}$.
 \vspace*{0.1cm}

\item If $X,Y$ are real normed vector spaces, we denote by $B(X,Y)$ the set (vector space)
of the linear, bounder operators from $X$ into $Y$.
 \vspace*{0.1cm}

 \item We denote by $\mathfrak{F}$ the Fourier transform on $L^2(\R^N)$, normalized in such a way
 that it is an \emph{isometry}; as a consequence, for every $f\in L^2(\R^N)\cap L^1(\R^N)$ we have
 $$\mathfrak{F}(f)(\xi) = \frac{1}{(2\pi)^{N/2}}\int_{\R^N}e^{-\imath\langle x,\xi\rangle}f(x)\,dx.$$
\end{itemize}

As anticipated in the Introduction,
in this `preliminary' section we collect several definitions and known results,
which will allow us to clearly state our main contribution (see Theorem
\ref{thm:Main} in Section \ref{mr}), and to
make the manuscript as self-contained as possible.

\subsection{The mixed operator $\LL = -\Delta+(-\Delta)^s$}
In order to clearly state the main result of this paper, we first need
to fix some notation and to pro\-perly define
what we mean by a \emph{solution to the Cauchy problem \eqref{eq:PbCauchy}}; due to the
\emph{mixed nature of $\LL$}, this will require some preliminaries.
\medskip

\noindent\textbf{1) The Fractional Laplacian.} Let $s\in (0,1)$ be fixed, and let
$u:\R^N\to\R$. The \emph{fractional La\-pla\-cian} (of order $s$) of $u$
at a point $x\in\R^N$ is defined as follows
\begin{equation} \label{eq:defDeltas}
\begin{split}
 (-\Delta)^s u(x) & = C_{N,s}\cdot \mathrm{P.V.}\int_{\R^N}\frac{u(x)-u(y)}{|x-y|^{N+2s}}\,dy
\\
& = C_{N,s}\cdot\lim_{\varepsilon\to 0^+}\int_{\{|x-y|\geq\varepsilon\}}\frac{u(x)-u(y)}{|x-y|^{N+2s}}\,dy,
\end{split}
\end{equation}
provided that the limit exists and is finite.
Here, $C_{N,s} > 0$ is a suitable normalization constant which plays
a role in the limit as $s\to 0^+$ or $s\to 1^-$, and is explicitly given by
$$C_{N,s} = \frac{2^{2s-1}{2s}\Gamma((N+2s)/2)}{\pi^{N/2}\Gamma(1-s)}.$$
As it is reasonable to expect, for $(-\Delta)^s u(x)$ to be well-defined one needs
to impose suitable \emph{growth conditions} on the function $u$, both when $|y|\to+\infty$ and
when $y\to x$. In this perspective we state the following
proposition (see \cite{KKP,Silv0} for a proof).

\begin{proposition} \label{prop:welldefDeltas}
Let $\Omega\subseteq\R^N$ be an open set. Then, the following facts hold.
\begin{itemize}
 \item[{i)}] If $0<s<1/2$ and $u\in C_{\mathrm{loc}}^{2s+\g}(\Omega)\cap {L}_s(\R^N)$
 for some $\g \in (0,1-2s)$, then
 $$\exists\,\,(-\Delta)^s u(x) = C_{N,s}\,\int_{\R^N}\frac{u(y)-u(x)}{|x-y|^{N+2s}}\,dy\quad
 \text{for all $x\in\Omega$}. $$
 \item[{ii)}] If $1/2<s<1$ and $u\in C^{1,2s-1+\g}_{\mathrm{loc}}(\Omega)\cap {L}_s(\R^N)$
 for some $\g\in (0,2-2s)$, then
 $$\exists\,\,(-\Delta)^s u(x) = -\frac{C_{N,s}}{2}\,\int_{\R^N}\frac{u(x+z)+u(x-z)-2u(x)}{|z|^{N+2s}}\,dy\quad
 \text{for all $x\in\Omega$}. $$
\end{itemize}
Moreover, in both cases \emph{i)\,-\,ii)} we have $(-\Delta)^su\in C(\Omega)$.
\end{proposition}
In the particular case when $\Omega = \R^N$
and $u\in\mathcal{S}\subseteq {L}_s(\R^N)$
(here and throughout, $\mathcal{S}$ denotes the usual Schwartz space of rapidly decreasing functions),
it is possible to provide an alternative expression of $(-\Delta)^s u$ (which is well-defined
on the whole of $\R^N$, see Proposition \ref{prop:welldefDeltas}) via the Fourier Transform
$\mathfrak{F}$; more precisely, we have the subsequent result.
\begin{proposition} \label{prop:DeltasFourier}
 Let $u\in\mathcal{S}\subseteq {L}_s(\R^N)$. Then,
 \begin{equation} \label{eq23}
  \exists\,\,(-\Delta)^s u(x) =
  \mathfrak F^{-1} \big( |\xi|^{2s}\mathfrak F u  \big)(x)\quad \text{for every
  $x\in\R^N$}.
 \end{equation}
\end{proposition}
 It should be noticed that, on account of \eqref{eq23}, it is immediate to recognize that
 Schwartz space $\mathcal{S}$ \emph{is not preserved} by the fractional Laplacian $(-\Delta)^s$
 (as $|\xi|^{2s}\mathfrak F u$ is not regular at $\xi = 0$), that is,
 one has $(-\Delta)^s(\mathcal{S})\not\subseteq\mathcal{S}$; however,
 we have the following characterization of the image
 $$\mathcal{S}_s = (-\Delta)^s(\mathcal{S}),$$
 which will be crucial to give the definition of \emph{solution of problem \eqref{eq:PbCauchy}}.
 \begin{proposition}[{See, e.g., \cite[Lem.\,1]{Stinga}}] \label{prop:DeltasS}
  Setting $\mathcal{S}_s = (-\Delta)^s(\mathcal{S})$, we have
  $$\mathcal{S}_s = \big\{\psi\in C^\infty(\R^N):\,
   \text{$(1 + |x|^{N+2s})D^\alpha\psi\in L^\infty(\R^N)$ for every
   $\alpha\in(\mathbb{N}\cup\{0\})^N$}\big\}.$$
 \end{proposition}
 Another consequence of the `representation formula'
 \eqref{eq23}, which plays a fundamental role in our argument
 (and, in general, in the analysis of the fractional Laplace operator $(-\Delta)^s$),
 is the possibility of realizing this operator as a \emph{densely defined, self-adjoint and 
 non-negative}
 operator on the Hilbert space $L^2(\R^N)$, whose associated heat semigroup admits a global heat kernel.
 Indeed, taking into account \eqref{eq23}, it is natural to define
  \begin{equation} \label{eq:realizationBs}
  \begin{gathered}
  \mathcal{B}_s:H^s(\R^N)\subseteq L^2(\R^N)\to L^2(\R^N),
  \qquad \mathcal{B}_s(u) = \mathfrak{F}^{-1}\big(|\xi|^{2s}\mathfrak{F}(u)\big) \\
  \text{where $H^s(\R^N) =  \big\{u\in L^2(\R^N):\,|\xi|^{2s}\mathfrak{F}(u)\in L^2(\R^N)\big\}$}.
  \end{gathered}
  \end{equation}
  Clearly, we have $\mathcal{S}\subseteq H^s(\R^N)$, and thus $\mathcal{B}_s$ is densely defined;
  moreover, by \eqref{eq23} one has
  $$\mathcal{B}_s(u) = (-\Delta)^s u\quad\text{for every $u\in\mathcal{S}\subseteq H^s(\R^N)$},$$
  and this shows that $\mathcal{B}_s$ is indeed a realization of $(-\Delta)^s$ on $L^2(\R^N)$.
  We then observe that, since the map $\mathfrak{F}$ is an isometry of $L^2(\R^N)$,
  for every $u,v\in H^s(\R^N)$ we get
  \begin{alignat*}{2}
   \mathrm{i)}\,\,&\langle \mathcal{B}_s(u),v\rangle_{L^2(\R^N)}
   && = \langle\mathfrak{F}(\mathcal{B}_s(u)),\mathfrak{F}(v)\rangle_{L^2(\R^N)}
   = \langle|\xi^{2s}\mathfrak{F}(u),\mathfrak{F}(v)\rangle_{L^2(\R^N)}
   \\
   & && = \langle\mathfrak{F}(u),|\xi|^{2s}\mathfrak{F}(v)\rangle_{L^2(\R^N)}
   = \langle u,\mathfrak{F}^{-1}(|\xi|^{2s}\mathfrak{F}(v))\rangle_{L^2(\R^N)} \\
   & && = \langle u, \mathcal{B}_s(v)\rangle_{L^2(\R^N)}; \\[0.1cm]
   \mathrm{ii)}\,\,&\langle \mathcal{B}_s(u),u\rangle_{L^2(\R^N)}
   && = \langle\mathfrak{F}(\mathcal{B}_s(u)),\mathfrak{F}(u)\rangle_{L^2(\R^N)}
   = \langle|\xi^{2s}\mathfrak{F}(u),\mathfrak{F}(u)\rangle_{L^2(\R^N)}
   \\
   & && = \langle|\xi|^s\mathfrak{F}(u),|\xi|^{s}\mathfrak{F}(u)\rangle_{L^2(\R^N)}\geq 0;
  \end{alignat*}
  and thus $\mathcal{B}_s$ is self-adjoint and non-negative. As a consequence of these facts,
  we are then entitled to apply \cite[Thm.\,4.9]{GrigoryanBook}, ensuring that
   the operator $-\mathcal{B}_s$ generates a strongly continuous semigroup
  on the Hilbert space $L^2(\R^N)$, say $\{T(t)\}_{t\geq 0}$. By this, we mean that
  \begin{enumerate}
   \item[P1)] for every fixed $t\geq 0$, we have $T(t)\in{B}(L^2(\R^N),L^2(\R^N))$;
   \item[P2)] $T(t+\tau) = T(t)\circ T(\tau)$ for every $t,\tau\geq 0$;
   \item[P3)] for every fixed $t\geq 0$ and $f\in L^2(\R^N)$, we have
   $$\lim_{\tau\to t}T(\tau)f = T(t)f\quad \text{in $L^2(\R^N)$};$$
   \item[P4)] for every fixed $t > 0$ and $f\in L^2(\R^N)$, we have
   $T(t)f\in H^s(\R^N)$ and
   $$\frac{d}{dt}\big(T(t)f\big) =
    \lim_{h\to 0}\frac{T(t+h)f-T(t)f}{h}= -\mathcal{B}_s\big(T(t)f\big)\,\,\text{in $L^2(\R^N)$}.$$
  \end{enumerate}
  This semigroup is called the
  \emph{heat semigroup of $-(-\Delta)^s$}, and it is denoted by $(e^{-t(-\Delta)^s})_{t\geq 0}$.
  \vspace*{0.1cm}

 We now observe that,
 starting from property P4) and exploiting the Fourier transform (together with the very definition
 of $\mathcal{B}_s$), it is easy to show that the operator
 $e^{-t(-\Delta)^s}$ (for every $t > 0$) is actually
 an \emph{integral operator on $L^2(\R^N)$ with a kernel of convolution type.}
 \vspace*{0.1cm}

 Indeed, let $f\in L^2(\R^N)$ be fixed, and let
 $$u:[0,+\infty)\to L^2(\R^N),\qquad
 u(t)(x) =  e^{-t(-\Delta)^s}f(x).$$
 Using property P4) and applying the Fourier transform, we see that
 \begin{alignat*}{2}
   \ast)&\,\,\mathfrak{F}(u'(t)) && =
   \mathfrak{F}\Big(x\mapsto\frac{d}{dt}\big(e^{-t(-\Delta)^s}f\big)(x)\Big)
   = -\mathfrak{F}\big(x\mapsto\mathcal{B}_s(e^{-t(-\Delta)^s}f)(x)\big) \\
   &\,&& = -|\xi|^{2s}\mathfrak{F}\big(x\mapsto e^{-t(-\Delta)^s}f(x)\big)
   = -|\xi|^{2s}\mathfrak{F}(u(t)), \\
   \ast)&\,\,\mathfrak{F}(u(0)) && = \mathfrak{F}\big(x\mapsto e^{-0\cdot(-\Delta)^s}f(x)\big)
   = \mathfrak{F}(f),
 \end{alignat*}
 which is a (formal) \emph{first-order, linear Cauchy problem} for $t\mapsto \mathfrak{F}(u(t))(\xi)$
 (for every fixed $\xi\in\R^N$); as a consequence,
 by formally solving this problem, we derive
 $$\mathfrak{F}(u(t))(\xi) = \mathfrak{F}(f)(\xi)e^{-t|\xi|^{2s}}\quad\text{for all $\xi\in\R^N,\,t\geq 0$}.$$
 Since we have expressed $\mathfrak{F}(u(t))$ as a \emph{product of two functions},
 by using the well-known properties of the Fourier transform we then conclude that
 \begin{equation} \label{eq:reprheatsemigroup}
  \begin{split}
   e^{-t(-\Delta)^s}f(x) & = u(t)(x) =
   \mathfrak{F}^{-1}\big(e^{-t|\xi|^{2s}}\cdot\mathfrak{F}(f)\big) \\
   & =
   (\mathfrak{h}_t^{(s)}*f)(x)
 = \int_{\R^N}\mathfrak{h}^{(s)}_t(x-y)f(y)\,dy,
  \end{split}
 \end{equation}
where, for every $z\in\R^N$ and $t > 0$, we have
 \begin{equation} \label{eq:fractionalHeatdef}
 \begin{gathered}
 \mathfrak{h}_t^{(s)}(z) = \frac{1}{(2\pi)^{N/2}}
 \mathfrak{F}^{-1}\big(e^{-t|\xi|^{2s}}\big)(z)
 = \frac{1}{(2\pi)^{N}}\int_{\R^N}e^{\imath\langle z,\xi\rangle-t|\xi|^{2s}}\,d\xi.
 \end{gathered}
 \end{equation}
 This function $(t,z)\mapsto \mathfrak{h}^{(s)}_t(z)$ is usually referred to as the
 \emph{heat kernel of $-(-\Delta)^s$}, and it satisfies
 the following properties (see, e.g., \cite{BarBassChenKass, BassLevin, 
 ChenKumagai, ChenKimKumagai, Sug} for a complete proof):
  \begin{enumerate}
   \item $\mathfrak{h}^{(s)}\in C^\infty(\R^+\times\R^N)$ and $\mathfrak{h}^{(s)} > 0$.
   \item For every $x\in\R^N$ and $t > 0$, we have
   $$\mathfrak{h}^{(s)}_t(x) = \mathfrak{h}^{(s)}_t(-x)\quad\text{and}\quad
    \mathfrak{h}^{(s)}_t(x) = \frac{1}{t^{N/(2s)}}\mathfrak{h}^{(s)}_1(t^{-N/(2s)}x).$$
    \item For every fixed $x\in\R^N$ and $t > 0$, we have
    $$\int_{\R^N}\mathfrak{h}^{(s)}_t(x)\,dy = 1.$$
    \item For every fixed $x\in\R^N$ and $t,\tau > 0$, we have
    $$\int_{\R^N}\mathfrak{h}^{(s)}_t(x-y)\mathfrak{h}^{(s)}_{\tau}(y)\,dy =
    \mathfrak{h}^{(s)}_{t+\tau}(x).$$
    \item There exists
   $C \geq 1$ such that
  \begin{equation} \label{eq:estimHeatFractional}
  \begin{gathered}
   C^{-1}\min\Big\{t^{-N/(2s)},\frac{t}{|x|^{N+2s}}\Big\}\leq \mathfrak{h}^{(s)}_t(x)
   \leq C\min\Big\{t^{-N/(2s)},\frac{t}{|x|^{N+2s}}\Big\} \\
   \text{for every $x\in\R^N$ and every $t > 0$}.
   \end{gathered}
  \end{equation}
  \end{enumerate}
\medskip

\noindent \textbf{2)\,\,The heat kernel of $\LL$.}
Now we have reviewed a few basic concepts on the fractional Laplace operator $(-\Delta)^s$,
we spend a few words concerning the heat semingroup and the associated global heat kernel of
the operator $-\LL = \Delta - (-\Delta)^s$
(we refer, e.g., to \cite{SV} for a thorough investigation on
this topic); this kernel will be used to introduce the notion of
\emph{mild solution}
to the Cauchy problem \eqref{eq:PbCauchy} (see Definition \ref{def:solWeakMild} below).

\vspace*{0.1cm}

Our starting point is the usual realization of the operator $-\Delta$
in $L^2(\R^N)$: deno\-ting by $H^2(\R^N)$ the classical Sobolev space $W^{2,2}(\R^N)$, it is
very well-known that
the operator
$$\mathcal{A}:H^2(\R^N)\to L^2(\R^N),\qquad \mathcal{A}(u) = \mathfrak{F}^{-1}\big(|\xi|^2\mathfrak{F}(u)\big),$$
satisfies the following properties:
\begin{itemize}
 \item[a)] $\mathcal{A}$ is a densely defined, positive and self-adjoint operator;
 \item[b)] $\mathcal{A}(u) = -\Delta u$ for every $u\in\mathcal{S}\subseteq H^2(\R^N)$
\end{itemize}
(actually, the above properties of $\mathcal{A}$ can be proved
by repeating \emph{verbatim} the computation carried out in the previous paragraph
with the `formal' choice $s = 1$, see also \cite[Sec.\,4.3]{Evans}).

On the other hand, by exploiting
the characterization of the Sobolev spaces $H^k(\R^N)$ (for $k\geq 1$) in terms
of $\mathfrak{F}$ (see, e.g., \cite[Sec.\,5.8.4]{Evans}), we have
$$H^2(\R^N) = \big\{u\in L^2(\R^N):\,|\xi|^{2}\mathfrak{F}(u)\in L^2(\R^N)\big\}\subseteq H^s(\R^N);$$
thus, taking into account \eqref{eq:realizationBs}, we can define
\begin{equation*}
 \mathcal{P}: H^2(\R^N)\to L^2(\R^N), \qquad \mathcal{P}(u) = \mathcal{A}(u)+\mathcal{B}_s(u) =
 \mathfrak{F}^{-1}\big(|\xi|^2\mathfrak{F}(u)\big)
  +\mathfrak{F}^{-1}\big(|\xi|^{2s}\mathfrak{F}(u)\big) .
\end{equation*}
Clearly, by combining the properties
of $\mathcal{B}_s$ (discussed in the previous paragraph)
with the pro\-perties of $\mathcal{A}$ recalled above, we immediately derive that
$\mathcal{P} = \mathcal{A}+\mathcal{B}_s$ is
a densely defined,  positive and self-adjoint operator which realizes $\LL$
on  $L^2(\R^N)$: indeed, we have
$$\mathcal{P}(u) = \LL u\quad\text{for every $u\in\mathcal{S}\subseteq H^2(\R^N)$}.$$
  We can then exploit once again \cite[Thm.\,4.9]{GrigoryanBook}, which ensures
  that also the operator $-\mathcal{P}$ generates a
  \emph{strongly continuous semigroup in the Hilbert space $L^2(\R^N)$}, which we denote by
  $$(e^{-t\LL})_{t\geq 0}$$
  (that is, the family $(e^{-t\LL})_{t\geq 0}$ satisfies
  the same properties P1)\,-\,to\,-\,P4) in the previous paragraph, with
  $-\mathcal{P}$ in place of $-\mathcal{B}_s$); this semigroup
  is called the \emph{heat semigroup of $-\LL$}.
  \vspace*{0.1cm}

  Now, by arguing exactly as in the previous paragraph, we see that the operator $e^{-t\LL}$ (for every
  fixed $t > 0$) is a integral operator with convolution-type kernel;
  more precisely, we have
  \begin{equation} \label{eq:reprSemigroupLL}
   \begin{gathered}
 e^{-t\LL}f(x) = (\mathfrak{p}_{t}*f)(x)
 = \int_{\R^N}\mathfrak{p}_{t}(x-y)f(y)\,dy\quad (\text{for every $f\in L^2(\R^N)$}),
 \end{gathered}
 \end{equation}
 where, for every $z\in\R^N$ and $t > 0$, we have
 \begin{equation} \label{eq:heatkernelLLI}
   \begin{gathered}
    \mathfrak{p}_t(z) = \frac{1}{(2\pi)^{N/2}}
 \mathfrak{F}^{-1}\big(e^{-t(|\xi|^{2}+|\xi|^{2s})}\big)(z)
 = \frac{1}{(2\pi)^{N}}\int_{\R^N}e^{\imath\langle z,\xi\rangle-t(|\xi|^2+|\xi|^{2s})}\,d\xi
 \end{gathered}
  \end{equation}
  (notice that $\mathfrak{F}(\mathcal{P}f) = (|\xi|^{2}+|\xi|^{2s})\mathfrak{F}(f)$);
  on the other hand, by exploiting \eqref{eq:fractionalHeatdef} (jointly with
  the e\-xplicit expression of $\mathfrak{F}^{-1}(e^{-t|\xi|^2})$ and the properties
  of the Fourier transform), we obtain
  \begin{equation*}
   \begin{split}
    \mathfrak{p}_t(z) & = \frac{1}{(2\pi)^{N/2}}
     \mathfrak{F}^{-1}\big(e^{-t|\xi|^{2}}\cdot e^{-t|\xi|^{2s}}\big)(z)
    \\
    & = \frac{1}{(2\pi)^{N/2}}\mathfrak{F}^{-1}\Big((2\pi)^{N/2}\mathfrak{F}(\mathfrak{g}_t)\cdot
    (2\pi)^{N/2}\mathfrak{F}(\mathfrak{h}_t^{(s)})\Big)(z) \\
    & = \mathfrak{F}^{-1}\Big((2\pi)^{N/2}\,\mathfrak{F}(\mathfrak{g}_t)\cdot
    \mathfrak{F}(\mathfrak{h}_t^{(s)})\Big)(z) \\
    & = (\mathfrak{g}_t*\mathfrak{h}_t^{(s)})(z),
   \end{split}
  \end{equation*}
  where $\mathfrak{g}_t(z)$ is the usual Gauss-Weierstrass heat kernel
  of $\Delta$, that is,
  $$\mathfrak{g}_t(z) = \frac{1}{(4\pi t)^{N/2}}e^{-|z|^2/(4t)}.$$
  Summing up, we conclude that
  \begin{equation} \label{eq:heatkernelLLII}
   \mathfrak{p}_t(z) = \frac{1}{(4\pi t)^{N/2}}
   \int_{\R^N}e^{-|z-\zeta|^2/(4t)}\mathfrak{h}^{(s)}(\zeta)\,d\zeta\quad (z\in\R^N,\,t > 0).
  \end{equation}
  This function $(t,z)\mapsto \mathfrak{p}_t(z)$ is referred to as
  the \emph{heat kernel of $-\LL$}, and it satisfies
 analogous properties to that of $\mathfrak{h}^{(s)}$;
 for a future reference, we collect these properties (which easily follow from
 the `explicit' expression of $\mathfrak{p}$ in \eqref{eq:heatkernelLLI}\,-\,\eqref{eq:heatkernelLLII})
 in the next theorem.
 \begin{theorem} \label{thm:HeatKerL}
 The heat kernel $\mathfrak{p}$ satisfies the following properties.
  \begin{enumerate}
      \item $\mathfrak{p}\in C^\infty(\R^+\times\R^N)$ and $\mathfrak{p} > 0$.
   \item For every $x\in\R^N$ and $t > 0$, we have
   $$\mathfrak{p}_t(x) = \mathfrak{p}_t(-x).$$
   \item For every fixed $x\in\R^N$ and $t > 0$, we have
    $$\int_{\R^N}\mathfrak{p}_t(x-y)\,dy = 1.$$
    \item For every fixed $x\in\R^N$ and $t,\tau > 0$, we have
    $$\int_{\R^N}\mathfrak{p}_t(x-y)\mathfrak{p}_\tau(y)\,dy =
    \mathfrak{p}_{t+\tau}(x).$$
    \end{enumerate}
   Moreover, by combining \eqref{eq:estimHeatFractional}
   with the `convolution-type' expression of $\mathfrak{p}_t$ in
   \eqref{eq:heatkernelLLII},
   we deduce the following
   upper estimate: there exists a constant $C > 0$ such that
  \begin{equation} \label{eq:upperEstimpt}
   0< \mathfrak{p}_t(x)\leq Ct^{-\frac{n}{2s}}\quad\text{for every $x\in\R^N,\,t > 0$}.
  \end{equation}
  \end{theorem}
  We finally point out that, starting from property P4) of the heat semigroup
  $(e^{-t\LL})_{t\geq 0}$, it is quite standard to prove that
  the \emph{unique solution} of the `abstract' $L^2$-Cauchy problem
  $$\begin{cases}
  \partial_t u = -\LL u+ f & \text{in $\R^N\times (0,+\infty)$} \\
  u(x,0) = u_0 & \text{for $x\in\R^N$}
  \end{cases}$$
  (for any fixed $f,u_0\in L^2(\R^N)$) is given by
  $$u(x,t) = e^{-t\LL}u_0(x)+ \int_0^t (e^{-(t-\tau)\LL}f)(x)\,d\tau;$$
  thus, by \eqref{eq:heatkernelLLII} we can rewrite this unique solution as follows
  \begin{equation} \label{eq:repruabstract}
   u(x,t) = \int_{\R^N}\mathfrak{p}_t(y)u_0(y)\,dy
  + \iint_{S_t}\mathfrak{p}_{t-\tau}(x-y)f(y)\,dy\,d\tau.
  \end{equation}
\begin{remark} \label{rem:convolutionpt}
 It is worth mentioning that the `convolution-type' formula \eqref{eq:heatkernelLLII} of $\mathfrak{p}$
 can be easily proved by taking into account the \emph{probabilistic interpretation of the operator $\LL$}.

 Indeed, since $\LL$ is the sum of the two operators $-\Delta$ and $(-\Delta)^s$,
 it is the infinitesimal generator of a stochastic process, say $(X_t)_{t\geq 0}$, which is the
 sum of two \emph{independent processes}, namely a Brownian motion $(W_t)_{t\geq 0}$
 and a pure jump L\'evy flight $(J_t)_{t\geq 0}$;
 thus, given any $t > 0$, we know that the law of the process $X_t$ (which is the
 function $\mathfrak{p}_t$) is the convolution of the laws of $W_t$ (the Gauss-Weierstrass heat kernel
 $\mathfrak{g}_t$)
 and of $J_t$ (the fractional heat kernel $\mathfrak{h}_t^{(s)}$).
\end{remark}
\begin{remark} \label{rem:formalcomputationhp}
 It is important to stress that the computations carried out in the previous paragra\-phs
 in order to obtain the `explicit' expressions of $\mathfrak{h}_t^{(s)}$
 and of $\mathfrak{p}_t$ in \eqref{eq:fractionalHeatdef}\,-\,\eqref{eq:heatkernelLLII},
 respectively, are actually \emph{formal computations}; however,
 \emph{starting from the mentioned expressions \eqref{eq:fractionalHeatdef}\,-\,\eqref{eq:heatkernelLLII}},
 one can prove \emph{a posteriori} that
 all the properties of $\mathfrak{h}^{(s)}$ and of $\mathfrak{p}$ hold.
\end{remark}
\medskip

\section{Existence and non-existence results}\label{mr}\setcounter{equation}{0}
\subsection{Very weak and mild solutions to problem \eqref{eq:PbCauchy}.}
Taking into account all the facts recalled so far, we
can now make precise the notion of \emph{solution to the Cauchy problem \eqref{eq:PbCauchy}}.
Actually, as is customary in the context of parabolic problems, we consider two different notions
of solutions, that is, \emph{very weak} and \emph{mild}.
\begin{definition} \label{def:solWeakMild}
 Let $u_0\in L^\infty(\R^N),\,u_0\geq 0$, and let $1\leq p<+\infty$.
 \begin{itemize}
  \item[1)] (\textbf{Very weak solution}) We say that a function $u:\overline{S}\to\R^+_0$
 is a \textbf{very weak solution} to problem \eqref{eq:PbCauchy} if the following properties hold:
 \begin{itemize}
  \item[a)$_1$] $u\in L^p_{\mathrm{loc}}(\overline{S})$;
  \item[b)$_1$] given any $T > 0$, we have $u\in L^\infty((0,T);L_s(\R^N))$;
  \item[c)$_1$] given any $\varphi\in\mathcal{T}_0$, we have
  \begin{equation} \label{eq:weakSoldef}
   \iint_{S}u(-\partial_t\varphi+\LL\varphi)\,dx\,dt
   -\int_{\R^N}u_0(x)\varphi(x,0)\,dx = \iint_{S}u^p\varphi\,dx\,dt.
  \end{equation}
 \end{itemize}
 \item[2)] (\textbf{Mild solution})
 We say that a function $u:\overline{S}\to\R^+_0$
 is a
 \textbf{mild solution} to problem \eqref{eq:PbCauchy} if the following properties hold:
 \begin{itemize}
  \item[a)$_2$] $u\in C(\overline{S})\cap L^\infty(S)$;
  \item[b)$_2$] for every $(x,t)\in S$, we have the identity
  \begin{equation} \label{eq:mildSoldef}
   u(x,t) = \int_{\R^N}\mathfrak{p}_t(x-y)u_0(y)\,dy
  + \iint_{S_t}\mathfrak{p}_{t-\tau}(x-y)u^p(y,\tau)\,dy\,d\tau.
  \end{equation}
 \end{itemize}
 \end{itemize}
\end{definition}
\begin{remark} \label{rem:defSol}
We list, for a future reference, some
remarks concerning Definition \ref{def:solWeakMild}.
\begin{enumerate}
 \item Taking into account Proposition \ref{prop:DeltasS}, it is easy
 to check that identity \eqref{eq:weakSoldef} \emph{is meaningful}, that is,
 for every fixed test function $\varphi\in \mathcal{T}_0$ we have
 \begin{itemize}
  \item[i)] $u(-\partial_t\varphi+\LL \varphi),\,u^p\varphi\in L^1(S)$;
  \item[ii)] $u_0\varphi(\cdot,0)\in L^1(\R^N)$
  \end{itemize}
 (provided that $u$ satisfies properties a)$_1$\,-\,to\,-\,c)$_1$).
 \vspace*{0.05cm}

  In fact, let $r,T > 0$ be such that $\varphi\equiv 0$ out of $B_r\times[0,T)$.
 First of all we observe that, since by
 property {a)}$_1$ one has $u\in L^p(B_r\times (0,T))$, we immediately get
 \begin{align*}
  \iint_{S}|u^p\varphi|\,dx\,dt
  & \leq \|\varphi\|_{L^\infty(S)}\iint_{B_r\times(0,T)}u^p\,dx\,dt < +\infty.
 \end{align*}
 On the other hand, recalling that $\varphi\in\mathcal{T}_0$,
 using Proposition \ref{prop:DeltasS}
 (and taking into account the explicit proof
 of this proposition given in \cite[Thm.\,9.4]{CRSTV}) we derive that
 $$|(-\Delta)^s(x\mapsto \varphi(x,t))|\leq \frac{c}{1+|x|^{N+2s}}\,\mathbf{1}_{[0,T)}(t)
 \quad\text{for every $(x,t)\in S$},$$
 for some constant $c > 0$ independent of $t$;
 as a consequence, since $-\partial_t\varphi-\Delta\varphi$ is (smooth and)
 supported in $B_r\times [0,T)$, and since $u\in L^\infty((0,T);L_s(\R^N))$, we obtain
 \begin{align*}
  & \int_{S}|u(-\partial_t\varphi+\LL\varphi|\,dx\,dt \\
  & \qquad \leq \int_{B_r\times(0,T)}u|\partial_t\varphi+\Delta\varphi|\,dx\,dt
  +  c\int_0^T\Big(\int_{\R^N}\frac{u}{1+|x|^{N+2s}}\,dx\Big)dt \\
  & \qquad \leq c\Big(\|u\|_{L^1(B_r\times(0,T))}+\int_0^T\|u(\cdot,t)\|_{1,s}\,dt\Big) \\
  & \qquad\leq c\big(\|u\|_{L^1(B_r\times(0,T))}+\|u\|_{\infty,(0,T),L_s(\R^N)}\big) <+\infty,
 \end{align*}
 where we have used the fact that $u\in L^p(B_r\times (0,T))\subset L^1(B_r\times (0,T))$, and
 $c > 0$ is a constant (possibly different from line to line) only depending on $\varphi$.

 Finally, since $u_0\in L^\infty(\R^N)$ and $\varphi(\cdot,0)\in C_0^\infty(\R^N)$, we immediately
 infer that
 $$u_0\varphi(\cdot,0)\in L^1(\R^N).$$

 \item  Owing to the properties of $\mathfrak{p}$
 in Theorem \ref{thm:HeatKerL}, it is easy to check that
 also identity \eqref{eq:mildSoldef} is \emph{is meaningful}
 (provided that $u$ satisfies properties a)$_2$\,-\,b)$_2$). In fact,
 since by assumption we have $u_0\in L^\infty(\R^N)$, for every $x\in\R^N$ we get
 $$0\leq \int_{\R^N}\mathfrak{p}_t(x-y)u_0(y)\,dy \leq \|u_0\|_{L^\infty(\R^N)}
 \int_{\R^N}\mathfrak{p}_t(x-y)\,dy
 = \|u_0\|_{L^\infty(\R^N)}<+\infty.
 $$
 Moreover, since by property a)$_2$ we also have $u\in L^\infty(S)$, for every $(x,t)\in S$ we get
 \begin{equation} \label{eq:estimsecondTermine}
 \begin{split}
  & 0\leq \iint_{S_t}\mathfrak{p}_{t-\tau}(x-y)u^p(y,\tau)\,dy\,d\tau
  \\
  &\qquad \leq \|u\|^p_{L^\infty(S)}\int_0^t\Big(\int_{\R^N}\mathfrak{p}_{t-\tau}(x-y)\,dy\Big)d\tau
  = \|u\|^p_{L^\infty(S)}t < +\infty.
 \end{split}
 \end{equation}
 We explicitly notice that the definition of mild solution
 comes from the representation of the unique solution of the $L^2$-Cauchy problem
 for $\LL$ discussed
 in the previous paragraph: indeed, our Cauchy problem \eqref{eq:PbCauchy} can be rewritten as
 $$\begin{cases}
  \partial_t u = -\LL u + f & \text{in $\R^N\times(0,+\infty)$} \\
  u(x,0) = u_0(x) & \text{for every $x\in\R^N$};
 \end{cases}$$
 where $f = u^p$; hence, by the `representation formula'
 \eqref{eq:repruabstract} we should have
 $$u(x,t) = \int_{\R^N}\mathfrak{p}_t(y)u_0(y)\,dy
  + \iint_{S_t}\mathfrak{p}_{t-\tau}(x-y)f(y)\,dy\,d\tau,$$
 which is precisely formula \eqref{eq:mildSoldef} (with $f = u^p$).
 \vspace*{0.1cm}

 \item In the particular case when $u_0\in L^\infty(\R^N)\cap L^2(\R^N)$,
 if $u\in C(\overline{S})\cap L^\infty(S)$
 is any mild so\-lution of the Cauchy problem \eqref{eq:PbCauchy} it is easy to recognize that
 $$u(x,0) = u_0(x)\quad\text{for every $x\in\R^N$}.$$
 Indeed, since $u_0\in L^2(\R^N)$, by exploiting property P3) of the heat semigroup $(e^{-t\LL})_{t\geq 0}$,
 together with the representation \eqref{eq:reprSemigroupLL} and estimate \eqref{eq:estimsecondTermine},
 we get
 \begin{align*}
  & \lim_{n\to +\infty}u(x,1/n) \\
  & \qquad = \lim_{n\to +\infty}\Big(\int_{\R^N}\mathfrak{p}_{1/n}(x-y)u_0(y)\,dy+
  \iint_{S_t}\mathfrak{p}_{1/n-\tau}(x-y)u^p(y,\tau)\,dy\,d\tau\Big) \\
  &\qquad = \lim_{n\to +\infty}\Big(\int_{\R^N}\mathfrak{p}_{1/n}(x-y)u_0(y)\,dy\Big) \\
  & \qquad= \lim_{n\to +\infty}(e^{-1/n\LL}u_0)(x) = u_0(x)\qquad\text{for a.e.\,$x\in \R^N$}
  \end{align*}
  (up to a sub-sequence, since $e^{-t\LL}u_0\to u_0$ as
   $t\to 0^+$ in $L^2(\R^N$)); thus, since $u\in C(\overline{S})$, we
   infer that $u(x,0) = u_0(x)$ for (a.e.) $x\in\R^N$.
   In particular,
  by modifying $u_0$ on a set of zero Lebesgue measure if needed, we conclude that
  $$u_0\in C(\R^N)\quad\text{and}\quad u(x,0) = u_0(x)\,\,\text{for every $x\in\R^N$}.$$

 \item
 Owing to the properties of the heat kernel $\mathfrak{p}$ in Theorem
 \ref{thm:HeatKerL}, and adapting the approach in
 the proof of
 \cite[Lem.\,2.1]{BPT}, 
 it is not difficult to recognize that
 \emph{any mild solution of problem \eqref{eq:PbCauchy} is also a
 very weak solution}.
\end{enumerate}
\end{remark}
\medskip

\subsection{The main result}
Now we have properly introduced the two types of solutions
for the Cauchy problem \eqref{eq:PbCauchy} we are interested in, we are
finally ready to state the main result of this paper.
\begin{theorem} \label{thm:Main}
 Let $u_0\in L^\infty(\R^N),\,u_0\geq 0$, and let $1< p<\infty$. We define
 $$\overline{p} = 1+\frac{2s}{N}.$$
 Then, the following facts hold.
 \begin{itemize}
  \item[1)] (\textbf{Non-existence}) If $1< p\leq \overline{p}$, \emph{there do not exist global in time
  very weak solutions} to the Cauchy problem
  \eqref{eq:PbCauchy} with $u_0\not \equiv 0$.
  \vspace*{0.1cm}

  \item[2)] (\textbf{Global existence}) If $p > \overline{p}$, there exist
  $\delta_0, \tau_0 > 0$ such that the Cauchy problem \eqref{eq:PbCauchy} possesses
  at least one global in time \emph{very weak solution}, provided that
\begin{equation}\label{eq:Stima_su_u0}
u_0(x)<\delta_0 \, \mathfrak{p}_{\tau_0}(x),\quad \textrm{for a.e. } x\in \mathbb{R}^{N}.
\end{equation}
 \end{itemize}
\end{theorem}
\begin{remark} \label{rem:existenceMild}
 As it will be clear from the proof
 of Theorem \ref{thm:Main}\,-\,2), the solution we are able to construct
 in the case $p > \overline{p}$ is actually a \emph{mild solution}
 to the Cauchy problem \eqref{eq:PbCauchy}.
\end{remark}

\section{Proof of Theorem \ref{thm:Main}}\label{proof}\setcounter{equation}{0}
In this section we provide the full proof of Theorem \ref{thm:Main}.
To ease the readability, we establish the two assertions 1) and 2) (\emph{non-existence} and
\emph{global existence}) separately.
\begin{proof}[Proof of Theorem \ref{thm:Main}\,-\,1)]
(\textbf{Non-existence}) Let $1< p \leq \overline{p}$ be fixed, and
suppose that there exists a very weak solution
 of the Cauchy problem \eqref{eq:PbCauchy} (in the sense of Definition \ref{def:solWeakMild},
 and for some initial condition $u_0\in L^\infty(\R^N),\,u_0\geq 0$). We then aim at proving that
 \begin{equation} \label{eq:toproveVeryWeak}
  \text{$u \equiv 0$ a.e.\,in $S$}.
 \end{equation}
 Once we know that \eqref{eq:toproveVeryWeak} holds, from \eqref{eq:weakSoldef} we infer that
 $$\int_{\R^N}u_0(x)\varphi(x,0)\,dx = \iint_{S}u(-\partial_t\varphi+\LL\varphi)\,dx\,dt
   -\iint_{S}u^p\varphi\,dx\,dt
   = 0\quad\forall\,\,\varphi\in\mathcal{T}_0,$$
   for which we derive that $u_0\equiv 0$ a.e.\,in $\R^N$.
   Hence, we turn to establish \eqref{eq:toproveVeryWeak}.
   To this end, it is convenient to distinguish the following two cases
$$
\text{(a)\,\,$1<p<\bar p$\qquad and\qquad (b) $p=\bar p$}.$$

\noindent \textbf{Case (a)}. To begin with, we choose two functions
 $\zeta\in C_0^\infty(\R^N),\,\psi\in C^\infty(\R^+_0)$
   such that
   \begin{itemize}
    \item[i)] $\zeta\equiv 1$ on $B_{1/2}$ and $\zeta\equiv 0$ out of $B_1$;
    \item[ii)] $\psi\equiv 1$ on $[0,1/2)$ and $\psi\equiv 0$ on $[1,+\infty)$;
    \item[iii)] $0\leq \zeta,\psi\leq 1$.
   \end{itemize}
   Then, we arbitrarily fix $r > 1$, and we define
   \begin{equation*}
   \begin{gathered}
    \xi_r(x) := \zeta^m\Big(\frac{x}{r}\Big),\qquad\phi_r(t) = \psi^m\Big(\frac{t}{r^{2s}}\Big) \\
    \text{where $m := \frac{2p}{p-1}$}.
   \end{gathered}
   \end{equation*}
   Since, obviously, we have $\varphi(x,t) = \xi_r(x)\phi_{r}(t)\in\mathcal{T}_0$,
   we are entitled to use this function $\varphi$ as a test function
   in \eqref{eq:weakSoldef}: recalling that (by assumption) $u_0\geq 0$ a.e. in $\R^N$, this gives
   \begin{equation} \label{eq:doveusarelastima}
   \begin{split}
     \iint_{S}u^p\varphi\,dx\,dt
     & =
     \iint_{S}u(-\partial_t\varphi+\LL\varphi)\,dx\,dt
   -\int_{\R^N}u_0(x)\xi_r(x)\,dx \\
   &  \leq \iint_{S}u(-\partial_t\varphi+\LL\varphi)\,dx\,dt \\
   &  = \iint_{S}u(-\xi_r\partial_t\phi_r-\phi_r\Delta\xi_r+\phi_r(-\Delta)^s\xi_r)\,dx\,dt.
    \end{split}
    \end{equation}
  We now turn to estimate the right-hand side of the above inequality.

  To this aim we first observe that
  \begin{equation} \label{eq:estimLocalPart}
  \begin{split}
   \mathrm{i)}\,\,&\Delta\xi_r = mr^{-2}\big[\zeta^{m-1}\Delta\zeta+(m-1)\zeta^{m-2}|\nabla\zeta|^2\big](x/r);\\
   \mathrm{ii)}\,\,&\partial_t\phi_r = mr^{-2s}\big[\psi^{m-1}\partial_t\psi\big](t/r^{2s}).
   \end{split}
  \end{equation}
  Moreover, since the function $G(z) = z^m$ is \emph{convex}, by
  \cite[Lem.\,3.2]{PV1} we have
  \begin{equation} \label{eq:estimNonLocalPart}
  \begin{split}
   (-\Delta)^s\xi_r & = (-\Delta)^s\big(G\circ (x\mapsto\zeta(x/r))\big)
   \leq m\zeta^{m-1}\Big(\frac{x}{r}\Big)(-\Delta)^s(x\mapsto\zeta(x/r)) \\
   & = \frac{m}{r^{2s}}\zeta^{m-1}(x/r)\big[(-\Delta)^s\zeta\big](x/r).
   \end{split}
  \end{equation}
  Thus, by combining \eqref{eq:estimLocalPart}\,-\,\eqref{eq:estimNonLocalPart}
  (and since $r > 1$), we obtain
  \begin{equation} \label{eq:estimfinaleRHS}
   \begin{split}
    & -\xi_r\partial_t\phi_r-\phi_r\Delta\xi_r+\phi_r(-\Delta)^s\xi_r
    \leq |\xi_r\partial_t\phi_r+\phi_r\Delta\xi_r|+\phi_r(-\Delta)^s\xi_r \\
    &
    \qquad \leq \mathbf{c}r^{-2s}\big(\zeta(x/r)\psi(t/r^{2s})\big)^{m-2}= \mathbf{c}r^{-2s}\varphi^{\frac{m-2}m}
     \\
   & \qquad = \mathbf{c}r^{-2s}\varphi^{1/p},
   \end{split}
  \end{equation}
  where we have also used the fact that $(-\Delta)^s\zeta\in\mathcal{S}_s$ (as $\zeta\in C_0^\infty(\R^N)$,
  see Proposition \ref{prop:DeltasS}).
  \vspace*{0.1cm}

  With estimate \eqref{eq:estimfinaleRHS} at hand, we can easily conclude the proof
  of
  \eqref{eq:toproveVeryWeak}: indeed, by combining the cited
  \eqref{eq:estimfinaleRHS} with the above estimate
  \eqref{eq:doveusarelastima}, and by using H\"older's inequality, we get
  \begin{align*}
   \iint_{S}u^p\varphi\,dx\,dt
   &  \leq
   \mathbf{c}r^{-2s}\iint_{S}u\varphi^{1/p}\,dx\,dt \\
   & (\text{since $\varphi$ is supported in $B_r\times[0,r^{2s})$}) \\
   & =
   \mathbf{c}r^{-2s}\int_0^{r^{2s}}\int_{B_r}u\varphi^{1/p}\,dx\,dt \\
   & \leq \mathbf{c}r^{-2s+(2s+N)\frac{p-1}{p}}\Big(\int_{S}u^p\varphi\,dx\,dt\Big)^{1/p};
  \end{align*}
  as a consequence, since $\varphi\equiv 1$ on $B_{r/2}\times[0,r^{2s}/2)$, we obtain
  \begin{equation} \label{eq:topasslimit}
   \int_0^{r^{2s}/2}\int_{B_{r/2}}u^p\,dx\,dt \leq \iint_{S}u^p\varphi\,dx\,dt
   \leq \mathbf{c}r^{N+2s-\frac{2sp}{p-1}}.
   \end{equation}
  On the other hand, since are assuming that $1<p<\overline{p}$, we have
  $$N+2s-\frac{2sp}{p-1} < 0;$$
  then, by letting $r\to +\infty$ in the above
  \eqref{eq:topasslimit}
  and by using the Monotone Convergence Theorem
   (recall that $r > 1$ was arbitrarily fixed, and $u\geq 0$ a.e.\,in $S$), we derive that
  $$\iint_{S}u^p\,dx\,dt = 0,$$
  from which we conclude that $u\equiv 0$ a.e.\,in $S$, as desired.
\medskip

\noindent\textbf{Case (b)}.
In this case, we use some ideas exploited in the proof of \cite[Thm.\,1]{FGS}.

First of all we observe that, if $p = \overline{p}$, we have
\begin{equation}\label{e5f}
\delta := -2s+(2s+N)\frac{p-1}{p}=0;
\end{equation}
thus, by arguing as in \textbf{Case (a)}, by \eqref{eq:topasslimit} and \eqref{e5f} we get
\[\int_0^{r^{2s}/2}\int_{B_{r/2}}u^p\,dx\,dt \leq \mathbf{c},\]
 for some constant $\mathbf{c}>0$ independent of $r$.
In particular,
by letting $r\to +\infty$ and by using the Monotone Convergence Theorem, we
can infer that $\displaystyle u\in L^p(\mathbb R^N\times (0, +\infty))$.
\vspace*{0.1cm}

We now define, for any $r>1, \beta>1$, the functions
   \begin{equation*}
   \begin{gathered}
    \xi_{r, \beta}(x) := \zeta^m\Big(\frac{x}{\beta r}\Big),\qquad\phi_r(t) = \psi^m\Big(\frac{t}{r^{2s}}\Big),
   \end{gathered}
   \end{equation*}
   where $\zeta,\psi$ and $m$ are as in the previous case.
   Clearly, $\varphi(x,t) = \xi_{r, \beta}(x)\phi_r(t)\in\mathcal{T}_0$,
   so we can use this function $\varphi$ as a test function
   in \eqref{eq:weakSoldef}: since $u_0\geq 0$, this gives
   \begin{equation} \label{e1f}
   \begin{split}
     \iint_{S}u^p\varphi\,dx\,dt
     & =
     \iint_{S}u(-\partial_t\varphi+\LL\varphi)\,dx\,dt
   -\int_{\R^N}u_0(x)\xi_{r, \beta}(x)\,dx \\
   &  \leq \iint_{S}u(-\partial_t\varphi+\LL\varphi)\,dx\,dt \\
   &  = \iint_{S}u(-\xi_{r, \beta}\partial_t\phi_r-\phi_r\Delta\xi_{r, \beta}+\phi_r(-\Delta)^s\xi_{r, \beta})\,dx\,dt.
    \end{split}
    \end{equation}
 Moreover, by arguing exactly as in \textbf{Case (a)}, we have the estimate
  \begin{equation} \label{e2f}
  \begin{split}
   \mathrm{i)}\,\,&|\Delta\xi_{r, \beta}(x)| \leq
   \mathbf{c}(\beta r)^{-2}\xi_{r,\beta}^{1/p}(x)\quad\text{for every $x\in\R^N$}; \\
   \mathrm{ii)}\,\,
   &(-\Delta)^s\xi_{r, \beta}(x)
   \leq \mathbf{c}(\beta r)^{-2s}\xi_{r,\beta}^{1/p}(x)
   \quad\text{for every $x\in\R^N$}; \\
   \mathrm{iii)}\,\,&|\partial_t\phi_r(t)| \leq \mathbf{c}r^{-2s}\phi_r^{1/p}(t)
   \cdot\mathbf{1}_{\{r^{2s}/2 <t <
   r^{2s}\}}(t)
   \quad\text{for every $t> 0$};
   \end{split}
  \end{equation}
By combining \eqref{e1f}\,-\,\eqref{e2f}, and by using H\"older's inequality, we
then get
  \begin{equation*}
 \begin{aligned}
   & \iint_{S}u^p\varphi\,dx\,dt
     \\
     & \qquad \leq
   \mathbf{c}r^{-2s}\iint_{S}u\phi_r^{1/p}\xi_{r,\beta}\cdot\mathbf{1}_{\{r^{2s}/2 <t <
   r^{2s}\}}(t)\,dx\,dt +
   \mathbf{c}(\beta r)^{-2s}\iint_{S}u\xi_{r, \beta}^{1/p}\phi_r(t)\,dx\,dt \\
   & \qquad
    \leq   \mathbf{c}r^\delta\beta^{\frac{Np}{p-1}} \left(\int_{\frac{r^{2s}}{2}}^{r^{2s}}
    \int_{B_{\beta r}} u^p\,dx dt \right)^{\frac 1p}+
     \mathbf{c} r^\delta\beta^{-2s+\frac{N(p-1)}p}\left( \int_0^{r^{2s}}\int_{B_{\beta r}}
   u^p\,dx dt \right)^{\frac 1p} \\
   & \qquad
    =   \mathbf{c}\beta^{\frac{Np}{p-1}} \left(\int_{\frac{r^{2s}}{2}}^{r^{2s}}
    \int_{B_{\beta r}} u^p\,dx dt \right)^{\frac 1p}+
     \mathbf{c}\beta^{-2s+\frac{N(p-1)}p}\left( \int_0^{r^{2s}}\int_{B_{\beta r}}
   u^p\,dx dt \right)^{\frac 1p},
  \end{aligned}
  \end{equation*}
  where we have used the fact that $\delta = 0$, see \eqref{e5f}.

   In particular, since $\varphi\equiv 1$ on $B_{\frac{r\beta}{2}}\times[0,r^{2s}/2)$, we obtain
  \begin{equation} \label{e8f}
  \begin{aligned}
   &\int_0^{r^{2s}/2}\int_{B_{\beta r/2}}u^p\,dx\,dt \leq \iint_{S}u^p\varphi\,dx\,dt\\
   &\qquad \leq \mathbf{c}
   \beta^{\frac{Np}{p-1}}
   \left( \int_{\frac{r^{2s}}{2}}^{r^{2s}}\int_{B_{\beta r}} u^p\,dx dt \right)^{\frac 1p}
     + \mathbf{c} \beta^{-2s+\frac{N(p-1)}p}\left( \iint_{S} u^p dx dt \right)^{\frac 1p}\,.
   \end{aligned}
   \end{equation}
   With \eqref{e8f} at hand, we can finally complete the proof of
   \eqref{eq:toproveVeryWeak} in this case.
   In fact, since we
   have already recognized that $u\in L^p(\R^N\times(0,+\infty))$,
   for any fixed $\beta\in (1, r)$ we have
\begin{equation}\label{e6f}
\begin{aligned}
&\lim_{r\to +\infty}\int_{\frac{r^{2s}}2}^{r^{2s}}\int_{B_{\beta r}} u^p\,dx dt \\
&\qquad=\lim_{r\to +\infty}\int_{0}^{r^{2s}}\int_{B_{\beta r}} u^p\,dx dt - \lim_{r\to +\infty}
\int_{0}^{\frac{r^{2s}}2}\int_{B_{\beta r}} u^p \,dx dt\\
& \qquad= \int_0^{+\infty} \int_{\mathbb R^N} u^p dx dt -  \int_0^{+\infty} \int_{\mathbb R^N} u^p  dx dt = 0\,.
\end{aligned}
\end{equation}
  On the other hand, since $p=\bar p$, we also have
  \begin{equation}\label{e9f}
  -2s+ \frac{N(p-1)}p = -\frac{2s(p-1)}p <0\,.
  \end{equation}
By virtue of \eqref{e6f}\,-\,\eqref{e9f},
letting $r\to +\infty$ and then $\beta\to +\infty$ in \eqref{e8f}, we obtain
 $$\iint_{S}u^p\,dx\,dt = 0,$$
  from which we deduce that $u\equiv 0$ a.e.\,in $S$, as desired.
  \end{proof}

\begin{proof}[Proof of Theorem \ref{thm:Main}\,-\,2)] (\textbf{Global existence})
We adapt to the present situation the line of arguments of the proof of \cite[Thm.\,1.1]{Pascucci}.
Let \eqref{eq:Stima_su_u0} be in force for some $\delta_0, \tau_0>0$ to be chosen later on, and let us introduce the following notation:
\begin{equation}\label{eq:Def_u_tilde_0}
\tilde{u}_{0}(x,t):= \int_{\mathbb{R}^{N}}\mathfrak{p}_{t}(x-y)u_{0}(y)\, dy
\end{equation}
\noindent and
\begin{equation}\label{eq:Def_Phi}
\Phi u (x,y) := \iint_{S_t}\mathfrak{p}_{t-\tau}(x-y)u^{p}(y,\tau)\, dy d\tau.
\end{equation}
Thanks to \eqref{eq:Stima_su_u0}, we have that
\begin{equation*}
\tilde{u}_{0}(x,t) \leq \delta_0 \, \int_{\mathbb{R}^{N}}\mathfrak{p}_{t}(x-y) \mathfrak{p}_{\tau_0}(y)\, dy = \delta_{0} \,\mathfrak{p}_{t+\tau_{0}}(x),
\end{equation*}
\noindent where in the last step we used Theorem \ref{thm:HeatKerL}-(4).\\
Exploiting \eqref{eq:Def_u_tilde_0} we now define the recursive sequence of functions $(\tilde{u}_n)_{n\in \mathbb{N}}$ as
\begin{equation}\label{eq:Def_tilde_u_n}
\tilde{u}_{n+1}(x,t) :=\tilde{u}_{0}(x,t) + \Phi \tilde{u}_{n}(x,t).
\end{equation}
By induction, we can prove that $(\tilde{u}_n)_{n\in \mathbb{N}}$ is monotone increasing. Indeed,
\begin{equation*}
\begin{aligned}
\tilde{u}_{1}(x,t) &= \tilde{u}_{0}(x,t) + \Phi \tilde{u}_{0}(x,t)\\
&= \tilde{u}_{0}(x,t) + \iint_{S_t}\mathfrak{p}_{t-\tau}(x-y)\tilde{u}_{0}^{p}(y,\tau)\, dy d\tau \geq \tilde{u}_{0}(x,t),
\end{aligned}
\end{equation*}
\noindent and, assuming $\tilde{u}_{n}\geq \tilde{u}_{n-1}$, and hence  $\tilde{u}^{p}_{n}\geq \tilde{u}^{p}_{n-1}$, we have
\begin{equation*}
\begin{aligned}
\tilde{u}_{n+1}(x,t) &= \tilde{u}_{0}(x,t) + \Phi \tilde{u}_{n}(x,t)= \tilde{u}_{0}(x,t) + \iint_{S_t}\mathfrak{p}_{t-\tau}(x-y)\tilde{u}_{n}^{p}(y,\tau)\, dy d\tau \\
&\geq \tilde{u}_{0}(x,t) + \iint_{S_t}\mathfrak{p}_{t-\tau}(x-y)\tilde{u}_{n-1}^{p}(y,\tau)\, dy d\tau = \tilde{u}_{n}(x,t).
\end{aligned}
\end{equation*}
In order to properly choose $\delta_0>0$, we further define the increasing (since $\delta_0>0$) sequence of real numbers $(\delta_n)_{n\in \mathbb{N}}$ as
\begin{equation*}
\delta_{n+1} := \delta_0 + \delta_{n}^{p},
\end{equation*}
\noindent If we choose $\delta_0 >0$ small enough, the sequence
$(\delta_n)_{n\in \mathbb{N}}$ is convergent, and therefore
there exists $M \in \mathbb{R}^{+}$ such that
\begin{equation}\label{eq:bound_Delta_n}
\delta_n \leq M \quad \textrm{for every } n\in \mathbb{N}.
\end{equation}
Our next goal is to choose $\tau_0>0$ such that
\begin{equation}\label{eq:claim_3.16_Pascucci}
\tilde{u}_{n}(x,t) \leq \delta_{n} \, \mathfrak{p}_{t+\tau_{0}}(x), \quad \textrm{for every } (x,t)\in S \textrm{ and for every } n\in \mathbb{N}.
\end{equation}
Before proceeding by induction, recalling that $p>\overline{p} = 1 + \tfrac{2s}{N}$ and thanks to both \eqref{eq:upperEstimpt} and Theorem \ref{thm:HeatKerL}-(4), we notice that
\begin{equation}\label{eq:conto_p_t+tau}
\begin{split}
& \iint_{S_t} \mathfrak{p}_{t-\tau}(x-y) \mathfrak{p}^{p}_{\tau+\tau_0}(y)\, dy d\tau \\
& \qquad
 \leq C^{p-1}\, \mathfrak{p}_{t+\tau_0}(x) \int_{0}^{+\infty} (\tau+\tau_0)^{-N(p-1)/(2s)}d\tau < \mathfrak{p}_{t+\tau_{0}}(x),
\end{split}
\end{equation}
\noindent provided that $\tau_0 > 0$ is large enough, namely
\begin{equation*}
\tau_0 > \left( C^{1-p} \left(\dfrac{N(p-1)}{2s}-1\right) \right)^{2s/(2s-N(p-1))}.
\end{equation*}
Let us now go through the induction procedure. Firstly,
\begin{equation*}
\begin{aligned}
\tilde{u}_{1}(x,t) &= \tilde{u}_{0}(x,t) + \iint_{S_t}\mathfrak{p}_{t-\tau}(x-y) u^p_{0}(y,\tau)\, dy d\tau \\
& \leq \delta_0 \, \mathfrak{p}_{t+\tau_0}(x) + \delta_{0}^{p}\, \iint_{S_t}\mathfrak{p}_{t-\tau}(x-y) \mathfrak{p}^{p}_{\tau+\tau_0}(y)\, dy d\tau\\
&\leq \left(\delta_{0} + \delta_{0}^{p}\right) \, \mathfrak{p}_{\t+\tau_{0}}(x) = \delta_1 \, \mathfrak{p}_{\t+\tau_{0}}(x).
\end{aligned}
\end{equation*}
Now, assuming that \eqref{eq:claim_3.16_Pascucci} holds for a certain $n\in \mathbb{N}$, it follows that
\begin{equation*}
\begin{aligned}
\tilde{u}_{n+1}(x,t) &= \tilde{u}_{0}(x,t) + \iint_{S_t}\mathfrak{p}_{t-\tau}(x-y) \tilde{u}^{p}_{n}(y,\tau) \, dy d\tau \\
&\leq \delta_{0} \, \mathfrak{p}_{t+\tau_{0}}(x) + \delta_{n}^{p} \iint_{S_t}\mathfrak{p}_{t-\tau}(x-y) \mathfrak{p}_{t+\tau_{0}}^{p}(y) \, dy d\tau \\
&\leq \left(\delta_0 + \delta_n^p\right) \mathfrak{p}_{\tau+\tau_0}(x) = \delta_{n+1}\, \mathfrak{p}_{\tau+\tau_0}(x),
\end{aligned}
\end{equation*}
\noindent where we exploited once again \eqref{eq:conto_p_t+tau}.\\
Combining \eqref{eq:claim_3.16_Pascucci} and \eqref{eq:bound_Delta_n}, we find that
\begin{equation*}
\tilde{u}_{n}(x,t) \leq M \, \mathfrak{p}_{t+\tau_{0}}(x), \quad \textrm{for every } (x,t)\in S \textrm{ and for every } n\in \mathbb{N}.
\end{equation*}
Let us now consider the function $u:= \sup \tilde{u}_{n}$.
By monotone convergence, $u$ satisfies \eqref{eq:mildSoldef}
and therefore $u$ is the desired global mild solution to \eqref{eq:PbCauchy}. In view of Remark \ref{rem:defSol}-(4), $u$ is also a global in time very weak solution to  \eqref{eq:PbCauchy}.
This closes the proof.
\end{proof}

\bigskip

\noindent{\bf{Acknowledgments.}}  The authors are members of the Gruppo Nazionale per l'Analisi Matematica, la Probabilit\`a e le loro Applicazioni (GNAMPA, Italy) of the Istituto Nazionale di Alta Matematica (INdAM, Italy). S.B and E.V. are partially supported by PRIN project 2022R537CS ``$NO^3$ - Nodal Optimization, NOnlinear elliptic equations, NOnlocal geometric problems, with a focus on regularity''.
F.P. is partially supported PRIN project 2022SLTHCE ``Geometric-analytic methods for PDEs and applications''. E.V. and F. P. are partially supported by Indam-GNAMPA projects 2024.

%and they are partically supported by GNAMPA. The third author acknowledges that this work is part of the PRIN project 2022 Geometric-analytic methods for PDEs and applications,
%ref. 2022SLTHCE, financially supported by the EU, in the framework of the "Next Generation EU initiative".


\begin{thebibliography}{999}

\bibitem{BB} C. Bandle, H. Brunner,
{\it Blowup in diffusion equations: A survey},
J. Comput. Appl. Math. {\bf 97}, (1998), 3--22.

\bibitem{BPT} C. Bandle, M.A. Pozio, A. Tesei,
{\it The Fujita exponent for the Cauchy problem in the hyperbolic space,}
J. Differential Equations {\bf 251}, (2011), 2143--2163.

\bibitem{BarBassChenKass}
 M.T. Barlow, R.F. Bass, Z-Q. Chen, M. Kassmann, 
 \emph{Non-local Dirichlet forms and symmetric jump processes}, 
 Trans. Amer. Math. Soc. \textbf{361}, (2009), 1963--1999.
 
\bibitem{BassLevin} R.F. Bass, D.A. Levin, 
\emph{Transition probabilities for symmetric jump processes}, 
Trans. Amer. Math. Soc. \textbf{354}, (2002), 2933--2953.

\bibitem{Biagi2} S. Biagi, S. Dipierro, E. Valdinoci, E. Vecchi,
{\it A Mixed local and nonlocal elliptic operators: regularity and maximum pronciples},
Comm. Partial Differential Equations {\bf 47} (2022), 585--629 .

\bibitem{Biagi3} S. Biagi, S. Dipierro, E. Valdinoci, E. Vecchi,
{\it  A Faber-Krahn inequality for mixed local and nonlocal operators}, J. Anal. Math. {\bf 150}(2), (2023), 405--448.

\bibitem{BiagiBN} S. Biagi, S. Dipierro, E. Valdinoci, E. Vecchi,
{\it  A Brezis-Nirenberg type result for mixed local and nonlocal operators }, submitted.
Available at: \url{https://arxiv.org/abs/2209.07502}

\bibitem{Biagi4} S. Biagi, D. Mugnai, E. Vecchi,
{\it A Brezis-Oswald approach for mixed local and nonlocal operators},
Commun. Contemp. Math. {\bf 26}(2), (2024), 2250057, 28 pp.

\bibitem{Biagi5} S. Biagi, E. Vecchi,
{\em Multiplicity of positive solutions for mixed local-nonlocal singular critical problems}, submitted.
Available at: \url{https://arxiv.org/abs/2308.09794}

\bibitem{Biagi6} S. Biagi, E. Vecchi,
{\em On the existence of a second positive solution to mixed local-nonlocal 
concave-convex critical problems}, submitted.
Available at: \url{https://arxiv.org/abs/2403.18424}

\bibitem{BonfE} M. Bonforte, J. Endal, {\it Nonlocal Nonlinear Diffusion Equations. Smoothing Effects,
Green Functions, and Functional Inequalities}, J. Funct. Anal. {\bf 284} (2023) 109831.

\bibitem{BFV}
M. Bonforte, A. Figalli, J.L. V\'{a}zquez,
{\em Sharp boundary behaviour of solutions to semilinear nonlocal elliptic equations},
Calc. Var.  Partial Differ. Equ. {\bf 57}, (2018), Article 57.

\bibitem{BonfII}  M. Bonforte, P. Ibarrondo, M. Ispizua, {\it The Cauchy-Dirichlet Problem for Singular Nonlocal Diffusions
on Bounded Domains}, DCDS-A {\bf 43} (2023) 1090-1142.

\bibitem{ChenKumagai}
Z-Q. Chen, T. Kumagai, 
\emph{Heat kernel estimates for stable-like processes on $d$-sets}, 
Stochastic Process. Appl. \textbf{108}, (2003), 27--62.
 
\bibitem{ChenKimKumagai} 
Z-Q. Chen, P. Kim, T. Kumagai, 
\emph{Global heat kernel estimates for symmetric jump processes}, 
Trans. Amer. Math. Soc. \textbf{363}, (2011), 5021--5055.

\bibitem{CRSTV}
\'O. Ciaurri, L. Roncal, P. R. Stinga, J. L. Torrea, J. L. Varona,
\emph{Nonlocal
discrete diffusion equations and the fractional discrete Laplacian,
regularity and applications},
Adv. Math. \textbf{330}, (2018), 688--738.


\bibitem{DeFMin}
 C. De Filippis, G. Mingione,
 {\em Gradient regularity in mixed local and nonlocal problems},
 Math. Ann. {\bf 388}, (2024), 261--328.

\bibitem{DL}  K. Deng, H.A. Levine,
{\em The role of critical exponents in blow-up theorems: the sequel},
J. Math. Anal. Appl. \textbf{243}, (2000), 85--126.

\bibitem{DPLV23} S. Dipierro, E. Proietti Lippi, E. Valdinoci, {\it (Non)local logistic equations with Neumann conditions}, Ann. Inst. H. Poincaré, Anal. Non Lin., (to appear).

\bibitem{DV21} S. Dipierro, E. Valdinoci, {\it
Description of an ecological niche for a mixed local/nonlocal dispersal: An evolution equation and a new Neumann condition arising from the superposition of Brownian and Lévy processes},
Physica A: Stat. Mech. Appl., {\bf 575} (2021) 126052.






\bibitem{Evans} L. C. Evans, "Partial Differential Equations", Second Edition, American Mathematical Society (2010).

\bibitem{FGS} A. Z. Fino, E. I. Galakhov, O. A. Salieva.  {\it Nonexistence of global weak solutions for evolution equations with fractional Laplacian}, Mathematical Notes  {\bf 108}  (2020), 877-883.


\bibitem{F} H. Fujita,
{\em On the blowing up of solutions of the Cauchy problem for $u_t=\Delta u+u^{1+\alpha}$},
 J. Fac. Sci. Univ. Tokyo Sect. I \textbf{13}, (1966), 109--124.

\bibitem{GarainKinnunen}
P. Garain, J. Kinnunen,
{\em On the regularity theory for mixed local and nonlocal quasilinear elliptic equations}, Trans. Amer. Math. Soc. {\bf 375}(8), (2022), 5393--5423.

\bibitem{GarainKinnunen2}
P. Garain, J. Kinnunen,
{\it On the regularity theory for mixed local and nonlocal quasilinear parabolic equations},
Ann. Sc. Norm. Super. Pisa Cl. Sci. {\bf 25}(1), (2024), 495--541.

 \bibitem{GarainKinnunen3}
 P. Garain, J. Kinnunen,
{\it Weak Harnack inequality for a mixed local and nonlocal parabolic equation},
J. Differential Equations {\bf 360}, (2023), 373--406.

\bibitem{GarainLindgren}
P. Garain, E. Lindgren,
{\it Higher H\"{o}lder regularity for mixed local and nonlocal degenerate elliptic equations},
Calc. Var. Partial Differential Equations {\bf 62}, 67, (2023).


\bibitem{GrigoryanBook}
A. Grigor'yan,
\emph{Heat Kernel and Analysis on Manifolds},
American Mathematical Society, Providence, RI; International Press, Boston, MA, 2009. xviii+482 pp.




\bibitem{MeGP1} G. Grillo, G. Meglioli, F. Punzo,
{\it Smoothing effects and infinite time blowup for reaction-diffusion equations: an approach via Sobolev and Poincar\'e inequalities},
J. Math. Pures Appl. {\bf 151 }, (2021), 99--131.

\bibitem{MeGP2} G. Grillo, G. Meglioli, F. Punzo, {\it Global existence of solutions and smoothing effects for classes of reaction-diffusion equations on manifolds,} J. Evol. Eq. {\bf 21}, (2021), 2339--2375.

%\bibitem{MeGPu3} G. Grillo, G. Meglioli, F. Punzo, {\it Blow-up and global existence for semilinear parabolic equations on infinite graphs}, preprint (2024).

\bibitem{MeGP3} G. Grillo, G. Meglioli, F. Punzo, {\it Global existence for reaction-diffusion evolution equations driven by the $ {\text{p}} $-Laplacian on manifolds}, Math. Eng. {\bf 5}, (2023), 1--38.

\bibitem{MeGP4} G. Grillo, G. Meglioli, F. Punzo,  {\it Blow-up versus global existence of solutions for reaction-diffusion equations on classes of Riemannian manifolds}, Ann. Mat. Pura Appl. {\bf 202}, (2023),  1255--1270.



\bibitem{GMPu} G. Grillo, M. Muratori, F. Punzo, {\it Blow-up and global existence for solutions to the porous medium equation with reaction and slowly decaying density}, J. Differential Equations {\bf 269}, (2020), 8918--8958.


\bibitem{Sun1} Q. Gu, Y. Sun, J. Xiao, F. Xu, {\it Global positive solution to a semi-linear parabolic equation with potential on Riemannian manifold},
 Calc. Var. Partial Differential Equations {\bf 59}, 170, (2020).

\bibitem{H} K. Hayakawa, {\em On nonexistence of global solutions of some semilinear parabolic differential equations}, Proc. Japan Acad. \textbf{49}, (1973), 503--505.

\bibitem{HKN} N. Hayashi, E. I. Kaikina and P. I. Naumkin, {\it  Asymptotics for fractional nonlinear heat equations}, J. London Math. Soc. {\bf 72}, (2005), 663--688.

\bibitem{IKK} K. Ishige, T. Kawakami and K. Kobayashi, {\it  Asymptotics for a nonlinear integral equation with a generalized heat kernel} J. Evol. Equ. {\bf 14}, (2014), 749-777.


\bibitem{KST} K. Kobayashi, T. Sirao, H. Tanaka, {\it On the growing up problem for semilinear heat equations}, J. Math. Soc. Japan {\bf 29}, (1977), 407--424.

\bibitem{KKP} J. Korvenp\"{a}\"{a}, T. Kuusi, Tuomo, G. Palatucci, {\it The obstacle problem for nonlinear integro-differential operators}, Calc. Var. Partial Differential Equations {\bf 55}, (2016), 1--30.

\bibitem{LaisterS} R. Laister, M. Sierzega {\it A blow-up dichotomy for semilinear fractional heat equations}, Math. Ann. {\bf  381}, (2021) 75--90.


\bibitem{Levine} H. A. Levine, {\it The role of critical exponents in blowup theorems}, SIAM Rev. {\bf 32}, (1990), 262--288.

\bibitem{MaMoPu} P. Mastrolia, D.D. Monticelli, F. Punzo, {\it Nonexistence of solutions to parabolic differential inequalities with a potential on Riemannian manifolds,} Math. Ann., {\bf 367} (2017), 929--963.

\bibitem{MeMoPu} G. Meglioli, D.D. Monticelli, F. Punzo, {\it Nonexistence of solutions to quasilinear parabolic equations with a potential in bounded domains,} Calc. Var. Partial Differential Equations {\bf 61}, 23, (2022).


\bibitem{Mitidieri2} E. Mitidieri, S.I. Pohozaev, {\it Towards a unified approach to nonexistence of solutions for a class of differential inequalities,} Milan J. Math. {\bf 72}, (2004), 129--162.

\bibitem{Pascucci} A. Pascucci,
{\it Semilinear equations on nilpotent Lie groups: global existence and blow-up of solutions},
Matematiche (Catania) {\bf 53}(2), (1998), 345--357.



\bibitem{Punzo} F. Punzo, {\it Blow-up of solutions to semilinear parabolic equations on Riemannian manifolds with
negative sectional curvature,} J. Math Anal. Appl. \bf 387 , \rm (2012), 815--827.

\bibitem{Pu22} F. Punzo, {\it Global solutions of semilinear parabolic equations with drift on Riemannian manifolds}, Discrete Contin. Dyn. Syst. {\bf 42}, (2022) 3733--3746.

\bibitem{PV1} F. Punzo, E. Valdinoci, {\it Uniqueness in weighted Lebesgue spaces for a class of fractional parabolic and elliptic equations}, J. Differential Equations {\bf 258}, (2015), 555--587.

\bibitem{PV2} F. Punzo, E. Valdinoci, {\it Prescribed conditions at infinity for fractional parabolic and elliptic equations with unbounded coefficients}, ESAIM: COCV {\bf 24}, (2018), 105--127.

\bibitem{QS}
 P. Quittner, P. Souplet,
 {\it Superlinear Parabolic Problems.
Blow-up, Global Existence and Steady States}, Birkh\"auser, Basel, 2007.

\bibitem{Silv0} L. Silvestre, {\it Regularity of the obstacle problem for a fractional power of the Laplace operator}, Ph.D. thesis, Univ. Texas at Austin (2006).

\bibitem{SV}
 R. Song. Z. Vondra\v{c}ek. 
 \emph{Parabolic Harnack Inequality for the Mixture of Brownian Motion and Stable Process}, 
 Tohoku Math. J. \textbf{59},  (2007), 1--19.

\bibitem{Sug} S. Sugitani, {\it
On nonexistence of global solutions for some nonlinear integral equations}
Osaka J. Math. {\bf 12}, (1975), 45--51.

\bibitem{Stinga}
P. R. Stinga,
{\it User's guide to the fractional Laplacian and the method of semigroups},
arXiv:1808.05159


\bibitem{WY2} Z. Wang, J. Yin, {\it Asymptotic behaviour of the lifespan of solutions for a semilinear heat equation in hyperbolic space}, Proc. Roy. Soc. Edinburgh Sect. A {\bf 146}, (2016), 1091--1114.

\bibitem{Wu} Y. Wu, {\it On nonexistence of global solutions for a semilinear heat equation on graphs}, Nonlinear Anal. {\bf 171}, (2018), 73--84.

\bibitem{Zhang} Q. S. Zhang, {\it Blow-up results for nonlinear parabolic equations on manifolds,} Duke Math. J. {\bf 97}(3), (1999), 515--539.








\end{thebibliography}
\end{document}